\documentclass[12pt]{amsart}
\usepackage{amsmath,amssymb,amsthm,mathrsfs,
amsfonts,array}
\newtheorem{theorem}{Theorem}[section]
\newtheorem{lemma}[theorem]{Lemma}
\newtheorem{definition}[theorem]{Definition}

\newtheorem{proposition}[theorem]{Proposition}
\theoremstyle{remark}
\newtheorem{remark}[theorem]{Remark}
\theoremstyle{hypothesis}
\newtheorem{hypothesis}[theorem]{Hypothesis}
\numberwithin{equation}{section}
\DeclareMathAlphabet\mbi{OML}{cmm}{b}{it}
\def\gm#1{\boldsymbol{#1}}
\def\mg#1{{\mbi{#1}}}
\def\eqdef{ \buildrel \triangle \over = }
\def\text#1{\mbox{#1}}

\def\exp#1{e^{#1}}

\newdimen\refbindent
\newcounter{parag}[subsection]
\newcounter{terme}[subsubsection]

\begin{document}
\title[Discrete-time multi-scale systems]
{Discrete-time multi-scale systems}
\author[D. Alpay]{Daniel Alpay}
\address{(DA) Department of mathematics\\ Ben-Gurion
University of the Negev\\ POB 653. Beer-Sheva
84105, Israel}
\email{dany@math.bgu.ac.il}
\author[M. Mboup]{Mamadou Mboup}
\address{(MM)  CReSTIC - Universit\'e
de Reims Champagne-Ardenne (and project ALIEN - Inria)\\ BP~1039, Moulin de la
Housse, 51687 Reims Cedex 2, France}
\email{Mamadou.Mboup@univ-reims.fr}
\date{}
\thanks{D. Alpay thanks the
Earl Katz family for endowing the chair which
supported his research. This research was
supported in part by the Israel Science
Foundation grant 1023/07}
\maketitle
\begin{abstract}
We introduce multi-scale filtering by the way of
certain double convolution systems. We prove
stability theorems for these systems and make
connections with function theory in the
poly-disc. Finally, we compare the framework
developed here with the white noise space
framework, within which a similar class of
double convolution systems has been defined
earlier.
\end{abstract}
\tableofcontents
\section{Introduction}
A wide class of causal discrete time-invariant
linear systems
can be given in terms of convolution in
the form
\begin{equation}
\label{1}
y_n=\sum_{m=0}^n h_{n-m}u_m,\quad n
=0,1,\ldots
\end{equation}
where $(h_n)$ is the impulse response and where
the input sequence $(u_m)$ and output sequence
$(y_m)$ belong to some sequences spaces. The
${Z}$ transform of the sequence $(h_n)$
\[
h(z)=\sum_{n=0}^\infty z^nh_n
\]
is called the transfer function of the system,
and there are deep relationships between
properties of $h$ and of the system. For instance
the system will be dissipative in the sense that
for all $\ell_2$ inputs $(u_n)$ it holds that
\[
\sum_{n=0}^\infty |y_n|^2\le \sum_{n=0}^\infty
|u_n|^2,
\]
if and only if $h$ is analytic and contractive
in the open unit disc ${\mathbb D}$. The function
$h$ is then called a {\sl Schur function}.
Similarly, the system will be $\ell_1-\ell_2$
bounded if for every $\ell_1$ entry, the output
is in $\ell_2$ and there is a $M>0$ independent
of the input such that
\[
\left(\sum_{n=0}^\infty |y_n|^2\right)^{1/2}\le
M\sum_{n=0}^\infty |u_n|.
\]
As is well known, the system is $\ell_1-\ell_2$
bounded if and only if $h$ belongs to the Hardy
space ${\mathbf H}_2({\mathbb D})$. Systems of
the form \eqref{1}, the Hardy space ${\mathbf
H}_2({\mathbb D})$ and Schur functions have been
generalized to a number of situations in the
theory of $N$-dimensional systems and beyond;
see for instance the works \cite{MR1839648},
\cite{MR2063749} \cite{btv}, \cite{MR2129642},
and the references therein.\\

Another generalization of systems of the form
\eqref{1} occurs when $h_n$ and $u_n$ are not
complex numbers, but belong to some space with a
product, say $\star$:
\begin{equation}
\label{system}
y_n=\sum_{m=0}^n h_{n-m}\star
u_m,\quad n=0,1,\ldots.
\end{equation}
Of special interest is the case where $\star$ is a
convolution  as, for instance in \cite{al_acap},
and the system is then called a
{\sl double convolution} system. Therein, the first
named author together
with David Levanony considered an example of
such a double convolution system, when both $h_n$
and $u_n$ are random variables, which belong to
the white noise space, or more generally to the
Kondratiev space. The product $h_{n-m}u_m$ in
\eqref{1} is then replaced by the Wick product.
The Wick product takes the form of a convolution
with respect to an appropriate basis, and we
have an example of a double convolution system.
Using the Hermite transform, one can define a
generalized transfer function, which is a
function analytic in $z$ and in a countable
number of other variables (these variables take
into account the randomness). The white noise space
setting is reviewed in the last section of
this paper, with purpose the comparison between the
present paper and \cite{al_acap}.\\

In the present work we study another type of double convolution
system, which arises in the theory of multi-scale systems. We use
the approach of the second named author presented at the {\sl
Mathematical Theory of Networks and Systems} conference in 2006
in Kyoto, see \cite{mboup_mtns}, to define the multi-scale
version of the systems \eqref{1}. Let
\[
\gm{\varphi}=\begin{pmatrix}a&b\\c&d\end{pmatrix}\in
SU(1,1)\quad\mbox{{\rm and  let}}\quad
\varphi(z)=\frac{az+b}{cz+d}
\]
be the corresponding automorphism of ${\mathbb D}$. Following
\cite{mboup_mtns}, consider the map
\[
T_{\gm{\varphi}}(f)(z)=\frac{1}{cz+d}f(\varphi(z))
\]
where $f$ is analytic in ${\mathbb D}$. Then
$T_{\gm{\varphi}}f$ is also analytic in the open unit disc.
Let
\[
f(z)=\sum_{n=0}^\infty f_nz^n\quad{\rm and}
\quad
(T_{\gm{\varphi}}f)(z)=\sum_{n=0}^\infty
f_{n, {\gm{\varphi}}}z^n
\]
be the Taylor expansions of $f$ and
$T_{\gm{\varphi}}f$ respectively. Let
\[
{\mathbb N}=\left\{1,2,\ldots\right\}\quad{\rm
and}\quad {\mathbb
N}_0=\left\{0,1,2,\ldots\right\}.
\]
The map which associates to the sequence $(f_n)_{n\in{\mathbb
N}_0}$ the sequence $(f_{n,{\gm{\varphi}}})_{{n\in\mathbb N}_0}$
is called the {\sl scaling operation} (see the precise definition
in Section \ref{sect2}). Consider now a subgroup $\Gamma$ of
$SU(1,1)$, which represents the scales we will use to study the
signals and systems.  One associates to the sequence
$(f_n)_{n\in{\mathbb N}_0}$ its {\sl scale transform}
$(f_n(\gamma))_{n\in{{\mathbb N}_0, \gamma \in \Gamma}}$, which
is a function of $n\in{\mathbb N}_0$ and $\gamma\in\Gamma$ and
where we have written $f_n(\gamma)$ for $f_{n, \gm{\gamma}}$. In
the case which we will consider, $\Gamma$ will be indexed by
${\mathbb Z }^p$, but the resulting
setting is quite different from classical $ND$ theory.\\

The scale transform is the starting point of our approach
(initiated in \cite{mboup_mtns}) to multi-scale analysis in
discrete time. In opposition to wavelets, we propose a transform
which has on the same level both the time and the scale aspect.
Let us now elaborate on the differences of our approach to
wavelets. Recall that in continuous time the wavelet transform of
a signal $f$ is defined by
\[
Wf(u,s)=\int_{\mathbb R}f(t)\frac{1}{\sqrt{s}}
\psi(\frac{t-u}{s})^*dt,
\]
where $s$ is the scale parameter and $\psi$
is the mother wavelet; see \cite[pp.
78-79]{Mallat-francais}. The discrete time
version of the transform is obtained, for discrete
scales $s = a^j$, by discretizing the above integral as in
\[
Wf[n, a^j]=\sum_k f_k \psi_{k-n, j}^*,
\]
where $\psi_{n,j} = \frac{1}{\sqrt{a^j}}
\psi\left(\frac{n}{a^j}\right)$ and $f_k = f(k)$, assuming a
sampling period normalized by 1 (see \cite[pp.
88-89]{Mallat-francais}). The decomposition provides an
appropriate mathematical tool for signal analysis. In particular,
it makes it possible to extract the components of a given
discrete-time signal at a given scale on a discrete grid.
However, the question of defining the scale shift operator
(dilation) for purely discrete-time signals is dodged somehow.
Another point of departure from our approach is that the wavelet
transform (either continuous- or discrete-time) has only one
convolution as compared to the double time and scale convolution
considered in the present work (see
equation \eqref{double} below).\\

We define linear systems as
expressions of the form \eqref{system}:
\begin{equation*}
 y_n=\sum_{m=0}^nh_{n-m} \star u_m,
\end{equation*}
where $\star$ denotes the convolution in
$\Gamma$. Thus,
\begin{equation}
\label{double}
y_n(\gamma)=\sum_{m=0}^n
\left(\sum_{\varphi\in\Gamma}
h_{n-m}(\gamma\circ\varphi^{-1})u_m(\varphi)
\right),\quad \gamma\in\Gamma.
\end{equation}
When $\Gamma$ is trivial we recover
\eqref{1}.\\

There are parallel and analogies between the
theory of linear stochastic systems presented in
\cite{al_acap} and the theory developed here.
These parallels are pointed out in the sequel,
and serve as guide and motivation for some of the
proofs in the present paper. Still, there are
some differences between the statements and the proofs
of the stability theorems in \cite{al_acap} and
the proofs given here. These differences are
pointed out in the text. We note that a general
theory of double convolution systems is in
preparation, \cite{alm}. We also note that some of the
results presented here have been announced in \cite{am_cras}.\\

We now turn to the outline of this paper. It consists of 10
sections besides the introduction. In Section \ref{sect2} we
review the approach to discrete multi-scale systems presented in
\cite{mboup_mtns}. In Section \ref{scalesyst} we define the
systems which we will study in the paper, and the related notion
of transfer function. Section \ref{trigomoment} is of a review
nature. We discuss the trigonometric moment problem and related
reproducing kernel Hilbert spaces of the type introduced by de
Branges and Rovnyak. In Section \ref{onegen} we consider the case
where the sub-group has one generator and is infinite. We use the
classical one dimensional moment problem to associate to the Haar
measure of the dual group a uniquely defined measure on the unit
circle. This allows us to use function theory on the disc and on
the bi-disc. In Section \ref{trigopoly} we review some deep
results of Mihai Putinar, see \cite{putinar_indiana}, on positive
polynomials on compact semi-algebraic sets and their use to solve
the trigonometric moment problem on the poly-disc. In Section
\ref{pg} we consider the case where $\Gamma$ is no more cyclic
but has a finite number of generators. Although the results in
Section \ref{onegen} are particular cases of the ones in Section
\ref{pg} we have chosen to present both for ease of exposition.
The next three sections consider stability results: BIBO
stability is considered in Section \ref{bibo}, dissipative
systems are studied in Section \ref{dissip} and Section
\ref{l1l2} is devoted to $\ell_1-\ell_2$ stability. In the last
section we review the white noise space setting, and  present a
table and some remarks, which point out the analogies between the
setting in \cite{al_acap} and the
present work.\\

{\bf Acknowledgements:} It is a pleasure to thank
Professor Mihai Putinar for explaining to us the
solution of the moment problem in the case of the
poly-disc.

\section{Scaling operator for discrete-time
signals} \label{sect2} We briefly summarize the approach to
multi-scale systems presented in \cite{mboup_mtns}. We first note
the following: If $F(s)$, $\Re{(s)} \geqslant 0$, denotes the
Laplace transform of a continuous-time signal $f(t)$, $t
\geqslant 0$, then, for any $\alpha = 1/\beta > 0$,
$\sqrt{\alpha} F(\alpha s)$ is the Laplace transform of $f(\beta
t)$. Therefore, time scaling has a similar form in the frequency
domain. As opposed to the continuous-time case, time scaling is
not clearly defined in the discrete-time setting. Nevertheless,
the preceding remark is the key step to define a scaling operator
for discrete-time signal.  Consider the M\"obius transformation
\begin{equation*}
G_{\theta}(s) = \frac{\exp{i\theta} - s}{\exp{-i\theta} + s},
\quad |\theta| < \frac{\pi}{2}
\end{equation*}
 which maps conformally the open right half-plane $\mathbb{C}_+$ onto
 the open unit disc. To recall our definition of the
 scaling operator
 (see \cite{mboup_mtns} and also \cite{AlpMbo_nds09}), we
 note that the scale shift in $\mathbb{C}_+$,
$$S_\alpha : s \mapsto S_\alpha(s) = \alpha s,
\quad \alpha > 0$$
translates in the unit disc, via $G_{\theta}(s)$, into the
hyperbolic transformation
\begin{equation}\label{gamma}
\gamma_{\{\alpha\}}(z) = (G_{\theta}\circ
S_\alpha\circ G_{\theta}^{-1})(z) =
\frac{(\exp{i\theta}+\alpha \exp{-i\theta})\ z +
(1 - \alpha)}{(1 - \alpha) z +
(\exp{-i\theta}+\alpha \exp{i\theta})}.
\end{equation}
Any such transformation  maps the open unit
disc (resp. the unit circle) into itself. Now, the most general
linear transformation  which maps the open unit disc (resp.
the unit circle) into itself has the form
\begin{equation}\label{jinner}
\gamma(z) = \frac{\gamma_1 z + \gamma_2}{\gamma_2^* z +
\gamma_1^*}, \quad |\gamma_1|^2 - |\gamma_2|^2 = 1.
  \end{equation}
If $|\Re(\gamma_1)| > 1$, then the transformation is hyperbolic
\cite{Ford} and it takes the form
\begin{equation}
\label{multiplier} \frac{\gamma(z) - \xi_1}{\gamma(z) - \xi_2} =
\alpha_\gamma \frac{z - \xi_1}{z - \xi_2}, \quad \alpha_\gamma > 0,
\end{equation}
where $\xi_1 = \frac{\sqrt{[\Re(\gamma_1)]^2-1} +i
  \Im(\gamma_1)}{\gamma_2^*} \eqdef \frac{\lambda_\gamma}{\gamma_2^*}$
and $\xi_2 = -\frac{\lambda_\gamma^*}{\gamma_2^*}$ are the two fixed
points. The constant $\alpha_\gamma$ is called
the multiplier  of the transformation, see \cite[p. 15]{Ford},
and is given by
$$\alpha_\gamma = \frac{\Re(\gamma_1) -
\sqrt{[\Re(\gamma_1)]^2-1}}{\Re(\gamma_1) +
\sqrt{[\Re(\gamma_1)]^2-1}}.$$ Noting that $|\xi_1|=|\xi_2|=1$,
one may rearrange \eqref{multiplier} to obtain
$$
\frac{\lambda_\gamma - \lambda_\gamma^*
\exp{i\xi_\gamma}\gamma(z)}{1 +
\exp{i\xi_\gamma}\gamma(z)} = \alpha_\gamma \frac{\lambda_\gamma -
\lambda_\gamma^* \exp{i\xi_\gamma}z}{1 + \exp{i\xi_\gamma}z},
$$
where  $\exp{i\xi_\gamma} = \frac{\lambda_\gamma}{\gamma_2}$.
Dividing both sides of this equality by $|\lambda_\gamma|$ and setting
$\exp{i\theta_\gamma} = \frac{\lambda_\gamma}{|\lambda_\gamma|}$,
we recover \eqref{gamma} up to a rotation:
\begin{equation}
\label{gmultiplier}
\exp{i\xi_\gamma} \gamma(z) = (G_{\theta_\gamma}\circ S_{\alpha_\gamma}\circ
 G_{\theta_\gamma}^{-1})(\exp{i\xi_\gamma}z).
\end{equation}
Any hyperbolic transformation $\gamma$ of the form \eqref{jinner} is
thus (conformally) equivalent to a scale shift $S_{\alpha_\gamma}$
in $\mathbb{C}_+$.

In the sequel, we will be interested in Abelian
subgroups of hyperbolic transformations, but the
remainder of this section deals with general linear
transformations
\[
\varphi(z)=\frac{az+b}{cz+d}, \text{ with } ad -
bc = 1,
\]
from the open unit disk onto itself. To each such
transformation, we associate (in bold letters)
\[
\gm{\varphi} =\begin{pmatrix}a&b\\ c&d\end{pmatrix} \in
SU(1,1),
\]
and we define
\[
(T_{\gm{\varphi}}f)(z)=\frac{1}{cz+d}f(\varphi(z)).
\]

\begin{lemma}
Let $\gm{\varphi}_1$ and $\gm{\varphi}_2$ belong to $SU(1,1)$. Then
\begin{equation}
\label{comp}
T_{\gm{\varphi}_2}\circ T_{\gm{\varphi}_1}=T_{\gm{\varphi}_1\gm{\varphi}_2},
\end{equation}
and in particular for every $\gm{\varphi}\in SU(1,1)$ it holds that
\begin{equation}
T_{\gm{\varphi}^{-1}}=(T_{\gm{\varphi}})^{-1}.
\label{inverse}
\end{equation}
\end{lemma}
{\bf Proof:} We have
\[
(T_{\gm{\varphi}_2}(T_{\gm{\varphi}_1}f))(z)= \frac{1}{c_2z+d_2}
\frac{1}{c_1\varphi_2(z)+d_1}
f(\varphi_1(\varphi_2(z))).
\]
But
\[
\begin{split}
\frac{1}{c_2z+d_2} \frac{1}{
c_1\varphi_2(z)+d_1}&= \frac{1}{c_2z+d_2}
\frac{1}{c_1\frac{a_2z+b_2}{c_2z +d_2}+d_1}\\
&=\frac{1}{(c_1a_2+d_1c_2)z+c_1b_2+d_1d_2},
\end{split}
\]
which ends the proof since the second row of
\[
\begin{pmatrix}a_1&b_1\\c_1&d_1\end{pmatrix}
\begin{pmatrix}a_2&b_2\\c_2&d_2\end{pmatrix}
\]
is equal to
\[
\begin{pmatrix}
c_1a_2+d_1c_2&c_1b_2+d_1d_2
\end{pmatrix}.
\]
\mbox{}\qed\mbox{}\\

Let $f$ be analytic in the open unit disc. Then
$T_{\gm{\varphi}}f $ is also analytic in the open unit disc. The
mapping $T_{\gm{\varphi}}$ induces a mapping from the space of
sequences coefficients of power series of functions analytic in a
neighborhood of the origin into itself: if $f(z)=\sum_{n=0}^\infty
z^nx_n$ is the power series expansion  at the origin of the
function $f$ analytic in ${\mathbb D}$, then
\[
(T_{\gm{\varphi}}f)(z)=\sum_{n=0}^\infty x_{n,\gm{\varphi}}z^n.
\]

In view of \eqref{inverse} we have:

\begin{proposition}
Let $\gm{\varphi}\in SU(1,1)$, and let $(a_n)_{n\in{\mathbb N}_0}$
be a sequence of complex numbers such that
$\limsup_{n\longrightarrow \infty}|a_{n+1}|^{1/(n+1)}\le 1$. Then
there exists a sequence of complex numbers $(b_n)_{n\in{\mathbb
N}_0}$ such that
$$\limsup_{n\longrightarrow
\infty}|b_{n+1}|^{1/{(n+1)}}\le 1,$$
and
\begin{equation*}
a_n=b_{n,\gm{\varphi}},\quad n=0,1,\ldots
\end{equation*}
\end{proposition}

{\bf Proof:} Let $a(z)=\sum_{n=0}^\infty
a_nz^n$. It suffices to define a series $b_{n,\gm{\varphi}}$
via the formula
\[
\sum_{n=0}^\infty b_{n,\gm{\varphi}} z^n=(T_{\gm{\varphi}^{-1}}a)(z),
\]
and use \eqref{inverse}.
\mbox{}\qed\mbox{}\\

\begin{theorem}
\label{unit}
The operator $T_{\gm{\varphi}}$ is unitary from
${\mathbf H}_2(\mathbb D)$ onto itself with norm equal to $1$. It
is also continuous from ${\mathbf H}_\infty(\mathbb D)$ into
itself with norm equal to $1/(|d|-|c|)$.
\end{theorem}

{\bf Proof:} First recall the formula
\begin{equation}
\label{rk_moebius}
\frac{1-\varphi(z)\varphi(w)^*}{1-zw^*}=
\frac{1}{(cz+d)(cw+d)^*}
\end{equation}
where $z,w$ are in the domain of definition of
$\varphi$. Furthermore, recall that a function
$f$ defined in ${\mathbb D}$ is analytic there
and belongs to ${\mathbf H}_2(\mathbb D)$, with
$\|f\|_{{\mathbf H}_2(\mathbb D)}\le 1$, if and
only if the kernel
\begin{equation}
\label{k_f}
\frac{1}{1-zw^*}-f(z)f(w)^*
\end{equation}
is positive in $\mathbb D$; see for instance
\cite[Theorem 2.6.6]{MR1839648}. We now compute
for
\[
\Delta(z,w)=
\frac{1}{1-zw^*}-(T_{\gm{\varphi}}f)(z)
(T_{\gm{\varphi}}f(w))^*
\]
for $z,w\in{\mathbb D}$. Using \eqref{rk_moebius}
we can write:
\[
\begin{split}
\Delta(z,w)&=\frac{1}{(1-\varphi(z)\varphi(w)^*)(cz+d)
(cw+d)^*}-\\
&\hspace{5mm}-\frac{1}{(cz+d)
(cw+d)^*}f(\varphi(z))f(\varphi(w))^*\\
&=\frac{1}{(cz+d)
(cw+d)^*}\times\\
&\hspace{5mm}\times
\left\{\frac{1}{1-\varphi(z)\varphi(w)^*}
-f(\varphi(z))f(\varphi(w))^*\right\}.
\end{split}
\]
The kernel
\[
\frac{1}{(1-\varphi(z)\varphi(w)^*)}
-f(\varphi(z)f(\varphi(w))^*
\]
is positive in ${\mathbb D}$ since the kernel
\eqref{k_f} is positive there. It follows that
$\Delta(z,w)$ is positive in the open unit disc,
and thus, by \cite[Theorem 2.6.6]{MR1839648}, the
function $T_{\gm{\varphi}}(f)$ belongs to ${\mathbf
H}_2(\mathbb D)$ and has norm less or equal to
$1$. To prove that the norm is indeed equal to
$1$ we use \eqref{comp} and \eqref{inverse},
which imply that
\[
1\le\|T_{\gm{\varphi}}\|\cdot\|T_{\gm{\varphi}^{-1}}\|,
\]
which, together with the fact that both
$T_{\gm{\varphi}}$ and $T_{\gm{\varphi}^{-1}}$
have norm less that $1$ implies that
$\|T_{\gm{\varphi}}\|=1$. We now show that
$T_{\gm{\varphi}}$ is unitary. Let $f\in{\mathbf
H}_2({\mathbb D})$. We have:
\[
\begin{split}
\|f\|_{{\mathbf H}_2({\mathbb D})}&=
\|T_{\gm{\varphi}^{-1}}T_{\gm{\varphi}}(f)
\|_{{\mathbf H}_2({\mathbb D})}\\
&\le\|T_{\gm{\varphi}}(f)\|_{{\mathbf H}_2({\mathbb D})}\\
&\le\|f\|_{{\mathbf H}_2({\mathbb D})},
\end{split}
\]
since both $T_{\gm{\varphi}}$ and
$T_{\gm{\varphi}^{-1}}$ are contractive. It
follows that $T_{\gm{\varphi}}$ is unitary.\\

The second claim is easily verified.
\mbox{}\qed\mbox{}\\

We therefore associate to the signal
$\mg{x}=\{x_n\}_{n\in{\mathbb N}_0}\in\ell_2$ the signal
$\{x_{n,\gm{\varphi}}\}_{\substack{\hspace{-.8cm}
n\in{\mathbb N}_0\\
\gm{\varphi}\in SU(1,1)}}$ indexed by ${\mathbb N}_0\times
SU(1,1)$. In the sequel, we will simplify the notation by writing
$x_n(\varphi)$ in place of $x_{n,\gm{\varphi}}$. For fixed $m$,
$\{x_m(\varphi)\}$, with  $\gm{\varphi} \in SU(1,1)$, represents
a scale signal, that is, the observation of the signal $\{x_n\}$
at time $m$ through the scales $\gm{\varphi} \in SU(1,1)$. In the
sequel, we will consider (by convention) the zooming as
corresponding to the ``positive'' scales.

\begin{definition}
\label{scale_causal} The scale-causal projection
of $\{x_m(\varphi)\}, \gm{\varphi} \in SU(1,1)$
is given by the restriction of $\{x_m(\varphi)\}$
to the scales $\gm{\varphi}$ for which the
multiplier is strictly less than one:
$\alpha_\varphi < 1$.
\end{definition}

\begin{definition}
\label{gamma+}
Given a discrete subgroup $\Gamma$ of $SU(1,1)$
we denote by $\Gamma_+$ the set of
transformations consisting of the identity and of
the scales $\gm{\varphi}$ for which the
multiplier is strictly less than one:
$\alpha_\varphi < 1$. The system \eqref{system}
will be scale-causal if the elements
$h_n\in\ell_2(\Gamma_+)$.
\end{definition}
Given $\gamma$ and $\varphi$ two elements of $\Gamma$, we will say that
$\gamma$ succeeds $\varphi$ and will note $\varphi \preccurlyeq \gamma$,
if  $\gamma\circ \varphi^{-1} \in \Gamma_+$ that is:
\begin{equation*}
\varphi\preccurlyeq\gamma\quad\iff
\alpha_{\gamma\circ\varphi^{-1}} \leqslant 1.
\end{equation*}
\begin{proposition}
The relation $\preccurlyeq$ defines a total order in $\Gamma$.
\end{proposition}

{\bf Proof:} Since we assume that $\Gamma$ is Abelian, all the
transformations must have the same fixed points. The parameters
$\xi_\gamma$ and $\theta_\gamma$ in \eqref{gmultiplier} are
therefore constant. The proof then follows upon noting that the
multiplier $\alpha_{\gamma\circ \varphi}$ is given by:
$\alpha_{\gamma\circ \varphi} = \alpha_{\gamma} \alpha_{\varphi}$.
\mbox{}\qed\mbox{}\\

With this order we  obtain a bijection
\begin{equation*}
\gamma\mapsto \varrho(\gamma)
\end{equation*}
between $\Gamma$ and ${\mathbb Z}$, and one can
identify $\ell_2(\Gamma)$ and $\ell_2({\mathbb
Z})$ and $\ell_2(\Gamma_+)$ and $\ell_2({\mathbb
N}_0)$.
\begin{remark}
Using the isomorphism we introduce the following
definition:
\begin{definition}
The function $u(\gamma)$ from $\Gamma_+$ into
${\mathbb C}$ has finite support if
\[
N(u) \eqdef
\max\left\{\varrho(\gamma)\,\,\mbox{such that}\,\,
u(\gamma)\not=0\right\}<\infty.
\]
The support of the function $u$ is the interval
$[0,N(u)]\subset{\mathbb N}_0$.
\end{definition}
The results of this section remain valid if we
replace the Hardy space ${\mathbf H}_2(\mathbb
D)$ by its vector-valued version ${\mathbf H}_2(
\mathbb D)\otimes {\mathcal H}$, where
${\mathcal H}$ is some Hilbert space. One can
define in particular the scaling of random
sequences when ${\mathcal H}$ is a probability
space.
\end{remark}

\begin{remark}
In \cite[equation (17)]{mboup_mtns} another kind
of systems are considered, with only one
convolution. The main issue there is the notion
of scale-invariance in a stronger form, which will not be
considered here; see also \cite{a_m} for related
work.
\end{remark}

Finally, we note that one could consider systems
non-causal with respect to $n$, that is of the
form
\[
y_n=\sum_{\mathbb Z} h_{n-m}\star u_m.
\]
Thus, there are really four possibilities for
the various stability theorems we present,
depending on whether we have time causality or
not, and scale-causality or not. In this paper
we only give part of all possible results.\\
\section{Discrete-scale invariant
systems and signals}
\label{scalesyst}
The scaling operators $T_{\gm{\varphi}}$ form a
group of operators from the Hardy space ${\mathbf
H}_2({\mathbb D})$ onto itself. From now on, we
discretize the scale axis and restrict
$\gm{\varphi}$ to a discrete subgroup $\Gamma$
of $SU(1,1)$.
 We will
take $\Gamma$ Abelian (cyclic) and consisting of
hyperbolic transformations, and we denote by
$\widehat{\Gamma}$ its dual group.
Recall that $\widehat{\Gamma}$ is formed by the set of functions
$$\sigma : \Gamma \to \mathbb{T} \text{ such that }
\sigma(\iota) = 1 \text{ and }\forall \ \gamma,
\varphi, \ \sigma(\gamma\circ \varphi) = \sigma(\gamma)\sigma(\varphi),$$
where $\iota$ stands for the identity transformation.
The elements of $\widehat{\Gamma}$ are called characters
of the group $\Gamma$ (see \cite{HewRoss}).
We denote by
$\widehat\mu$ the Haar measure of
$\widehat{\Gamma}$, which is compact by the
Pontryagin duality \cite{HewRoss}. We recall the
definition of the Fourier transform on $\Gamma$
and of its inverse:
\begin{equation*}
\begin{split}
\widehat{x}(\sigma)&=\sum_{\gamma\in\Gamma}
x(\gamma)\sigma(\gamma)^*,\\
x(\gamma)&=\int_{\widehat{\Gamma}}
\widehat{x}(\sigma)\sigma(\gamma)d\widehat\mu(\sigma).
\end{split}
\end{equation*}
The Haar measure $d\widehat{\mu}$ is normalized so that
Plancherel's theorem holds:
\[
\|f\|^2_{\ell_2(\Gamma)} \eqdef \sum_{\gamma\in \Gamma }|f(\gamma)|^2 =
\int_{\widehat{\Gamma}}|\widehat{f}(\sigma)|^2
d\widehat{\mu}(\sigma) \eqdef
\|\widehat{f}\|^2_{\mathbf{L}_2(d\widehat{\mu})}.
\]
See \cite[Theorem 8.4.2 p. 123]{deitmar}.
\begin{definition}
A signal will be a sequence $\{u_n(\cdot)\}_{n\in{\mathbb N}_0}$
of elements of $\ell_2(\Gamma)$, and such that the condition
\begin{equation}
\label{sup} \sup_{n=0,1,\ldots}
\|u_n(\cdot)\|_{\ell_2(\Gamma)}<\infty
\end{equation}
holds.\\
A scale-causal signal will be a sequence
$\{u_n(\cdot)\}_{n\in{\mathbb N}_0}$ of elements of
$\ell_2(\Gamma_+)$, and such that the condition
\begin{equation}
\label{sup12}
\sup_{n=0,1,\ldots}
\|u_n(\cdot)\|_{\ell_2(\Gamma_+)}<\infty
\end{equation}
holds.\\

\end{definition}

In the sequel, we will impose the following stronger norm
constrains on a signal, besides \eqref{sup} or \eqref{sup12},
namely:
\begin{equation}
\label{l2}
\sum_{n=0}^\infty
\|u_n(\cdot)\|^2_{\ell_2(\Gamma)}<\infty,
\end{equation}
or
\begin{equation}
\label{l1}
\sum_{n=0}^\infty \|u_n(\cdot)\|_{\ell_2(\Gamma)}<\infty,
\end{equation}
and similarly for scale-causal signals.\\

We note the following: a dissipative filter cannot be effective
at all scales. At some stage, details cannot be seen. These
intuitive facts are made more precise in the following
proposition.
\begin{proposition}\mbox{}\\
$(1)$ Assume that the supports of the $u_n$  are uniformly
bounded. Then, \eqref{l2} is in force.\\
$(2)$ Assume that the support of $u_n$ is infinite for all $n$.
Then, the sum on the left side of \eqref{l2} diverges.
\end{proposition}
{\bf Proof:} Let $N$ be such that the support of all the
functions $\gamma\mapsto u_n(\gamma)$ is inside $[0,N]$. Then,
\[
\begin{split}
\sum_{n=0}^\infty\|u_n(\cdot)\|^2_{\ell_2(\Gamma)}&=\sum_{n=0}^\infty
\sum_{\varrho(\gamma)=0}^N|u_n(\gamma)|^2\\
&=\sum_{\varrho(\gamma)=0}^N\sum_{n=0}^\infty|u_{n}(\gamma)|^2\\
&\le N\|u_n\|_{\ell_2({\mathbb N}_0)},
\end{split}
\]
since (see Theorem \ref{unit}) the maps $T_{\gamma}$ are unitary
from ${\mathbf H}_2({\mathbb D})$ onto itself. \\

The second claim is proved similarly.\mbox{}\qed\mbox{}\\

An example of $(u_n)$ satisfying Condition $(1)$ of the preceding
proposition has been presented in the paper \cite{mboup_mtns},
where the corresponding group $\Gamma$ is Fuchsian. This
was used therein, to define the scale unit-pulse signal. A
similar condition was also considered by P. Yuditskii
\cite{Yuditskii} in the description of the direct integral
of spaces of character-automorphic functions.

\begin{definition}
An impulse response (resp. a scale-causal impulse response) will
be a sequence $\{h_n(\cdot)\}_{n\in{\mathbb N}_0}$ of elements of
$\ell_2(\Gamma)$ (resp. of $\ell_2(\Gamma_+)$) such that for
every $n\in {\mathbb N}_0$, the multiplication operator
\begin{equation}
\label{T_h}
\mathcal{M}_{h_n}:\quad u\mapsto h_n\star
u,\quad n=0,1,\ldots
\end{equation}
is bounded from $\ell_2(\Gamma)$ into itself
(resp. from $\ell_2(\Gamma_+)$ into itself) and
such that
\begin{equation}
\label{sup1} \sup_{n=0,1,\ldots}
\|h_n\|_{\ell_2(\Gamma)}<\infty \quad(resp.\quad
\sup_{n=0,1,\ldots}
\|h_n\|_{\ell_2(\Gamma_+)}<\infty).
\end{equation}
\label{def:mult}
\end{definition}

The systems that we consider here are
defined by \eqref{system}, that
is, by the double convolution \eqref{double}, that we recall below:
\begin{equation}
y_n(\gamma)=\sum_{m\in{\mathbb Z}} \left(\sum_{\delta\in\Gamma}
h_{n-m}(\gamma\circ\delta^{-1})x_m(\delta) \right).
\label{eq:conv}
\end{equation}
In view of \eqref{sup1} the series
\begin{equation}
\label{conv2}
H(z,\sigma)=\sum_{n=0}^\infty
z^n\widehat{h}_n(\sigma)
\end{equation}
converges in the ${\mathbf L}_2(d\widehat{\mu})$
norm for every $z\in{\mathbb D}$. Taking the
Fourier transform (with respect to $\Gamma$) of both
sides of \eqref{double} we obtain
\[
\widehat{y}_n(\sigma)=\sum_{m=0}^n
\widehat{h}_{n-m}(\sigma)
\widehat{x}_m(\sigma),
\]
where the equality is in the ${\mathbf L}_2(d\widehat{\mu})$
sense. Taking now the $Z$ transform we get
\begin{equation}
\label{tf} {Y}(z,\sigma)={ H}(z,\sigma)U(z,\sigma),
\end{equation}
where
\begin{equation}
\label{conv} {Y}(z,\sigma)=\sum_{n=0}^\infty z^n
\widehat{y}_n(\sigma)\quad{\rm and}\quad
U(z,\sigma)= \sum_{n=0}^\infty z^n
\widehat{u}_n(\sigma),
\end{equation}
and where, for every $z\in{\mathbb D}$ the
equality in \eqref{tf} is $\widehat{\mu}$-a.e.\\

The function $H(z, \sigma)$ can be seen as the transfer function
of the discrete-time scale-invariant system. Formula \eqref{conv2}
suggests to define and study hierarchies of transfer functions,
for which the functions $\widehat{h_n}$ depend on $\sigma$ in some
pre-assigned way (for instance, when they are polynomials in
$\sigma)$, or when the function $H(z,\sigma)$ is a rational
function of $z$ or of $\sigma$. In the next two sections, under
the hypothesis that the sub-group $\Gamma$ has a finite number,
say $p$,  of generators, we will associate to the system
\eqref{system} an analytic function of $p+1$ variables, which we
will call the {\sl generalized transfer function of the system}.

\section{The trigonometric moment problem}
\label{trigomoment}
%
We first gather some
well known facts on the trigonometric moment
problem in form of a theorem.

\begin{theorem}
Given an infinite sequence $\ldots,
t_{-1},t_0,t_1,\ldots$ of complex numbers such
that
\[
t_{-n}=t_n^*,\quad n=0,1,\ldots,
\]
there exists a positive measure $d\nu$ on
$[0,2\pi)$ such that
\[
t_n=\int_0^{2\pi}e^{-in\theta}d\nu(\theta),\quad n\in{\mathbb
N}_0,
\]
if and only if all the Toeplitz matrices
\[
\mathcal{T}_N=(t_{n-m})_{n,m=0,\ldots , N}
\]
are non-negative.
\label{moment_pb}
\end{theorem}

See for instance \cite[Theorem 2.7 p. 66]{knud}.
The measure is then unique (when normalized). We
also recall that the sequence $(t_n)$ and the
measure $d\nu$ are related by
\[
t_0+2\sum_{n=1}^\infty t_nz^n=
\int_0^{2\pi}\frac{e^{i\theta}+z}{e^{i \theta}-z}d \nu(\theta),
\]
and thus the function
\[
\Phi(z)=t_0+2\sum_{n=1}^\infty
t_nz^n=\int_0^{2\pi}\frac{e^{i\theta}+z}{e^{i
\theta}-z}d \nu(\theta)
\]
is analytic and has a positive real part in the
open unit disc. Using Stieltjes inverse formula,
one can recover $\nu$ from $\Phi$ via the
formula
\[
\lim_{r\rightarrow 1}\int_a^b{\rm
Re}\{\Phi(re^{i\theta})\}d\theta= \nu(b) - \nu(a_-),
\]
where we assume that $\nu$ is right continuous. We also recall
that the function
\[
\frac{\Phi(z)+\Phi(w)^*}{2(1-zw^*)}=
\int_0^{2\pi}\frac{d\nu(\theta)}{(e^{i\theta}-z)
(e^{i\theta}-w)^*}
\]
is positive for $z,w\in{\mathbb
C}\setminus{\mathbb T}$. We denote by ${\mathcal
L}_+(\Phi)$ the associated reproducing kernel
Hilbert space when $z$ and $w$ are restricted to
the open unit disc. The following result has
first been proved by de Branges and  Shulman;
see \cite{dbs}. In the statement, ${\mathbf
H}_2(d\nu)$ denotes the closed linear span in
${\mathbf L}_2(d\nu)$ of the functions $z^m$ for
$m\ge 0$.

\begin{theorem}
The space ${\mathcal L}_+(\Phi)$ consists of the
functions of the form
\[
\widetilde{h}(z)=\int_{0}^{2\pi}
\frac{h(e^{it})e^{it}}{
e^{it}-z}d\nu(t), \quad h\in{\mathbf H}_2(d\nu),
\]
with norm
\[
\|\widetilde{h}\|_{{\mathcal L}_+(\Phi)}=\|h\|_{{\mathbf
H}_2(d\nu)}.
\]
\end{theorem}

We now recall some results on the structure of
the space ${\mathbf H}_2(d\nu)$.

\begin{theorem}
\label{210709}
Assume that ${\mathbf
H}_2(d\nu)\not= {\mathbf L}_2(d\nu)$. Then
${\mathbf H}_2(d\nu)$ is a reproducing kernel
Hilbert space, and its reproducing kernel is of
the form
\[
\frac{A(z)A(w)^*-B(z)B(w)^*}{1-zw^*},
\]
where $A(z)$ and $B(z)$ are functions analytic
off the unit circle.
\end{theorem}

{\bf Proof:} We assume that ${\mathbf H}_2(d\nu)\not= {\mathbf
L}_2(d\nu)$, and let $h_0\in{\mathbf L}_2(d\nu)\ominus {\mathbf
H}_2(d\nu)$. Let $\alpha\in{\mathbb C}\setminus {\mathbb T}$ be
such that
\[
\int_{0}^{2\pi}\frac{h_0(\theta)d\nu(\theta)}{e^{i\theta}-\alpha}
\not =0.
\]
Let $p$ be a polynomial; then $R_\alpha
p\in {\mathbf H}_2(d\nu)$, where
\[
(R_\alpha p)(z)=\frac{p(z)-p(\alpha)}{z-\alpha}.
\]
Then
\[
\langle R_\alpha p\, ,\, h_0\rangle_{{\mathbf H}_2(d\nu)}=0,
\]
and therefore we obtain
\[
p(\alpha)= \frac{\displaystyle\int_0^{2\pi} \frac{p(e^{i\theta})}
{e^{i\theta}-\alpha}d\theta}{\displaystyle\int_{0}^{2\pi}
\frac{h_0(\theta)d\nu(\theta)}{e^{i\theta}-\alpha}}.
\]
Therefore the map $p\mapsto p(\alpha)$ is
continuous on the polynomials, and extends to a
continuous map to ${\mathbf H}_2(d\nu)$. Therefore
 ${\mathbf H}_2(d\nu)$ is a reproducing kernel Hilbert
subspace of ${\mathbf L}_2(d\nu)$. The proof is then
finished by using
\cite[Theorem 3.1 p. 600]{ad1}.
\mbox{}\qed\mbox{}\\
\section{The case of  one generator}
\label{onegen} In this section, we consider the
case of a cyclic group $\Gamma$, generated by a
hyperbolic transformation $\gm{\gamma}_0 \in
SU(1,1)$. Any transformation in $\Gamma$ is thus
of the form $\gamma_0^m \eqdef
\underbrace{\gamma_0 \circ
\ldots \circ \gamma_0}_{m \textrm{ times}}$, $m \in \mathbb{Z}$.
\begin{theorem}
There exists a positive measure $d\nu(\theta)$ on
$[0,2\pi)$ such that
\begin{equation}
\label{moment} \int_{\widehat{\Gamma}}
\sigma(\gamma_0^m)d\widehat{\mu}(\sigma)
=\int_{0}^{2\pi} e^{im\theta}d\nu(\theta),\quad
m\in{\mathbb Z}.
\end{equation}
\end{theorem}
{\bf Proof:} We use Theorem \ref{moment_pb}. Let
\[
t_m=\int_{\widehat{\Gamma}}\sigma(\gamma_0^m)^*
d\widehat{\mu}(\sigma),\quad m\in{\mathbb Z}.
\]
Since $|\sigma(\gamma_0)|=1$ for all $\sigma\in
\widehat{\Gamma}$ we have
\[
t_{\ell-m}=
\int_{\widehat{\Gamma}}\sigma(\gamma_0^\ell)^*
\sigma(\gamma_0^m)d\widehat{\mu}(\sigma)=\langle
\sigma(\gamma_0^m),
\sigma(\gamma_0^\ell)\rangle_{{\mathbf
L}_2(d\widehat{\mu})},\]
and therefore all the Toeplitz matrices
\[
\mathcal{T}_N=(t_{\ell-m})_{\ell,m=0,\ldots N}
\]
are non-negative. It follows from Theorem
\ref{moment_pb} that there exists a uniquely
defined measure $d\nu$ such that
\[
t_m=\int_{0}^{2\pi}e^{-im\theta}d\nu(\theta),
\quad m=0,1,2,\ldots,
\]
and hence we obtain \eqref{moment}.

\mbox{}\qed\mbox{}\\

\begin{remark}
The proof of the previous theorem formalizes
the intuitive idea that
one can make the ``change of variable''
\[
\sigma(\gamma_0)=e^{i\theta(\sigma)}.
\]
\end{remark}

\begin{theorem}
The linear map ${\mathbf I}$ which to
$\sigma(\gamma_0^m)$ associates the function
$z^m$:
\begin{equation}
\label{hermite1} {\mathbf I}(\sigma(\gamma_0^m))=
z^m,\quad m\in{\mathbb Z},
\end{equation}
is an isomorphism from ${\mathbf L}_2(d\widehat{\mu})$
into ${\mathbf L}_2(d\nu)$.
\end{theorem}

{\bf Proof:} For a function $f$ of the form
\begin{equation}
\label{daumesnil} f(\sigma)=\sum_{-N}^M
c_n\sigma(\gamma_0^n)\quad {\rm where}\quad N,M\in{\mathbb
N}_0\quad{\rm and} \quad c_n\in{\mathbb C},
\end{equation}
we have
\[
\begin{split}
\|f\|_{{\mathbf L}_2(d\widehat{\mu})}^2&=
\sum_{n,m=-N,\ldots, M}c_nc_m^*t_{m-n}\\
&=\sum_{n,m=-N,\ldots, M}c_nc_m^*\int_0^{2\pi}
e^{-i(m-n)\theta}d\nu(\theta)
\\
&=\int_0^{2\pi}|\sum_{n=-N}^Mc_ne^{in\theta}|^2d\nu(\theta)\\
&= \|{\mathbf I}(f)\|_{{\mathbf L}_2(d\nu)}^2.
\end{split}
\]
The result follows by continuity since such $f$ are dense in
${\mathbf L}_2(d\widehat{\mu})$. To verify this last claim we
note the following: By Plancherel's theorem, the map from
$\ell_2(\Gamma)$ onto ${\mathbf L}_2(d\widehat{\mu})$ which to
the sequence which consists only of zeros, except the $n$-th
element which is equal to $1$, associates the function
$\sigma(\gamma_0)^n$, extends to a unitary map.
\mbox{}\qed\mbox{}\\

We will be interested in particular in the
positive powers of $\gm{\gamma}_0$, which
correspond to zooming (we consider that the
multiplier of $\gm{\gamma}_0$, \textit{i.e.} the
associated scale $\alpha_{\gamma_0}$, is less
than 1).

\begin{definition}
We denote by ${\mathbf H}_2(d\widehat{\mu})$ the
closure in ${\mathbf L}_2(d\widehat{\mu})$ of
the functions $\sigma(\gamma_0)^n$,
$n=0,1,2,\ldots$. Similarly, we denote by ${\mathbf
H}_2(d\nu)$ the closure in ${\mathbf
L}_2(d\nu)$ of the functions $z^n$,
$n=0,1,2,\ldots$.
\end{definition}

Note that it may happen that \( {\mathbf
L}_2(d\widehat{\mu})={\mathbf
H}_2(d\widehat{\mu})\).\\

Following \cite{al_acap} we introduce the next definition.

\begin{definition}
The map ${\mathbf I}$ will be called the Hermite
transform.
\end{definition}

Recall that $\widehat{\Gamma}$ is compact and
therefore
\begin{equation}
\label{eq:inclusion} {\mathbf L}_2(d\widehat{\mu})\subset{\mathbf
L}_1(d\widehat{\mu}).
\end{equation}
In general the product of two elements $f$ and $g$ in
 ${\mathbf
L}_2(d\widehat{\mu})$ does not belong to  ${\mathbf
L}_2(d\widehat{\mu})$, and one cannot define ${\mathbf I}(fg)$,
let alone compare it with the product ${\mathbf I}(f){\mathbf
I}(g)$. On the other hand, we will need in the sequel only  the
case where at least one of the elements in the product $fg$
defines a bounded multiplication operator from ${\mathbf
L}_2(d\widehat{\mu})$ into itself; see Definition \ref{def:mult}
and the proof of Theorem \ref{tm:scd} for instance. This is
exploited in the next theorem.
\begin{theorem}
\label{hermite} Let $f\in{\mathbf L}_2(d\widehat{\mu})$ such that
the operator of multiplication by $f$ defines a bounded operator
from ${\mathbf L}_2(d\widehat{\mu})$ into itself. Then for every
$g$ in ${\mathbf L}_2(d\widehat{\mu})$ it holds that:
\begin{equation}
\label{ankara1}
{\mathbf I}(f g)={\mathbf I}(f){\mathbf I}(g).
\end{equation}
\end{theorem}
{\bf Proof:} We note that the multiplicative property
\eqref{ankara1} holds for $f$ and $g$ of the form
\eqref{daumesnil}. To prove the theorem we
 first assume that $g$ is of the form
\eqref{daumesnil}, and consider a sequence $(p_n)$ of elements of
the form \eqref{daumesnil}, converging to $f$ in the ${\mathbf
L}_2(d\widehat{\mu})$ norm. The function $g$ is in particular
bounded, and so $fg\in{\mathbf L}_2(d\widehat{\mu})$, and we can
write:
\[
\|fg-p_ng\|_{{\mathbf L}_2(d\widehat{\mu})}\le K
\|f-p_n\|_{{\mathbf L}_2(d\widehat{\mu})},
\]
where $K>0$ is such that $|g|\le K$. Thus, $fp_n$ tends in
${\mathbf
L}_2(d\widehat{\mu})$- norm to $fg$.\\

The function ${\mathbf I}(g)$ is bounded, and ${\mathbf
I}(f)\in{\mathbf L}_2(d\nu)$. Therefore:
\[
\begin{split}
\| {\mathbf I}(fg)-{\mathbf I}(f){\mathbf I}(g) \|_{{\mathbf
L}_2(d\nu)}&\le
 \|{\mathbf
I}(fg)-{\mathbf I}(p_ng)\|_{{\mathbf
L}_2(d\nu)}+\\
&\hspace{5mm}+\|{\mathbf I}(p_ng)-{\mathbf
I}(f){\mathbf I}(g)\|_{{\mathbf L}_2(d\nu)}\\
&=\|fg-p_ng\|_{{\mathbf
L}_2(d\widehat{\mu})}+\\
&\hspace{5mm}+\|{\mathbf I}(p_n){\mathbf I}(g)-{\mathbf
I}(f){\mathbf I}(g)\|_{{\mathbf L}_2(d\nu)}\\
&\le\|fg-p_ng\|_{{\mathbf
L}_2(d\widehat{\mu})}+\\
&\hspace{5mm}+K_1\|p_n-g\|_{{\mathbf L}_2(d\widehat{\mu})},
\end{split}
\]
where $K_1>0$ is such that $|{\mathbf I}(g)|\le K_1$, and where
we have used that
\[
{\mathbf I}(p_ng)={\mathbf I}(p_n){\mathbf I}(g)
\]
since both $p_n$ and $g$ are of the form \eqref{daumesnil}.
Hence, we
obtain \eqref{ankara1} for $f$ and $g$ as asserted.\\

Let now $f,g\in{\mathbf L}_2(d\widehat{\mu})$ be such that
$fg\in{\mathbf L}_2(d\widehat{\mu})$. Then, ${\mathbf I}(fg)$ is
well defined. Let $(q_n)$ be a sequence of elements of the form
\eqref{daumesnil}, converging to $g$ in the ${\mathbf
L}_2(d\widehat{\mu})$ norm. Then, by the preceding argument,
\[
{\mathbf I}(fq_n)={\mathbf I}(f){\mathbf I}(q_n),\quad\forall
n\in{\mathbb N}.
\]
In view of this equation and of the inclusion \eqref{eq:inclusion}
we can write:
\[
\begin{split}
\|{\mathbf I}(fg)-{\mathbf I}(f){\mathbf I}(g) \|_{{\mathbf
L}_1(d\nu)}&\le \| {\mathbf I}(fg)-{\mathbf I}(f){\mathbf I}(q_n)
\|_{{\mathbf L}_1(d\nu)}+\\
&\hspace{5mm}+ \| {\mathbf I}(f){\mathbf I}(q_n)-{\mathbf
I}(f){\mathbf I}(g) \|_{{\mathbf L}_1(d\nu)}\\
&= \| {\mathbf I}(f(g-q_n))
\|_{{\mathbf L}_1(d\nu)}+\\
&\hspace{5mm}+ \| {\mathbf I}(f){\mathbf I}(q_n-g) \|_{{\mathbf
L}_1(d\nu)}.
\end{split}
\]
By Cauchy-Schwartz inequality and by the isometry property of
${\mathbf I}$ on ${\mathbf L}_2(d\widehat{\nu})$, we have that
\[
\begin{split}
\|{\mathbf I}(f(g-q_n))\|_{{\mathbf
L}_1(d\nu)}&\le\sqrt{\nu((0,2\pi])}\cdot \|{\mathbf
I}(f(g-q_n))\|_{{\mathbf L}_2(d\nu)} \\&=\sqrt{\nu((0,2\pi])}\cdot
\|f(g-q_n)\|_{{\mathbf L}_2(d\widehat{\mu})}.
\end{split}
\]
Since multiplication by $f$ is assumed to define a bounded
operator from ${\mathbf L}_2(d\widehat{\mu})$ into itself, there
exists a constant $C>0$ such that
\[
\|f(g-q_n)\|_{{\mathbf L}_2(d\widehat{\mu})}\le C
\|g-q_n\|_{{\mathbf L}_2(d\widehat{\mu})} \rightarrow 0\quad
{\rm as}\quad n\rightarrow\infty.
\]
Using once more the Cauchy-Schwartz inequality we have:
\[
\begin{split}
\|{\mathbf I}(f){\mathbf I}(q_n-g)\|_{{\mathbf
L}_1(d{\nu})}&\le\|{\mathbf I}(f)\|_{{\mathbf
L}_2(d{\nu})}\cdot\|{\mathbf I}(q_n)-{\mathbf I}(g)\|_{{\mathbf
L}_2({\nu})}\\
&= \|f\|_{{\mathbf L}_2(d\widehat{\mu})}\cdot\|q_n-g\|_{{\mathbf
L}_2(d\widehat{\mu})}\rightarrow 0\quad {\rm as}\quad
n\rightarrow\infty.
\end{split}
\]
The claim follows.
\mbox{}\qed\mbox{}\\

At this stage we need a change of notation; since two (and, in
the next section, $p+1$) complex variables appear, we denote by
$z_1$ (and by $z_1,\ldots , z_p$ in the following section) the
variables related to the Hermite transform, and keep $z$ for the
$Z$-transform variable (this notation differs from the one in
\cite{al_acap}, where the $Z$-transform variable is denoted by
$\zeta$).

\begin{definition}
\label{dn:gtf}
The function
\begin{equation}
\mathscr H(z,z_1)=\sum_{n=0}^\infty z^n{\mathbf
I}(\widehat{h_n})(z_1)
\end{equation}
is called the {\sl generalized transfer function}
of the system.
\end{definition}

Taking the Hermite transform on both sides of \eqref{tf}, or,
equivalently, taking the $Z$ transform and the Hermite transform
on both sides of \eqref{eq:conv}, we obtain
\[
\mathscr Y(z,z_1)=\mathscr H(z,z_1)\mathscr U(z,z_1),
\]
where $\mathscr U(z,z_1)=\sum_{n=0}^\infty z^n{\mathbf
I}(\widehat{u_n})(z_1)$, and similarly for $\mathscr Y(z,z_1)$.
The function $\mathscr H$ is analytic in a neighborhood of
$(0,0)\in{\mathbb C}^2$. It is of interest to relate the
properties of $\mathscr H$ and of the system. This is done in
Sections $9$ and $10$ of the paper. We first study, in the next
section, the case where $\Gamma$ has a finite number of
generators.

\section{The trigonometric moment problem in the
poly-disc}
\label{trigopoly}
In \cite{putinar_indiana} a solution is given to
the $K$-moment problem when $K$ is a compact
semi-algebraic set. The material is quite deep,
and cannot be easily summarized in a short
overview here. The purpose of this section is to
serve as a guide to the reader to the topic. The
starting point is a semi-algebraic subset of
${\mathbb R}^n$, defined by the positivity of
$\kappa$ polynomials
\[
K=\left\{x\in{\mathbb R}^n\,;\, p_j(x)\ge 0\, ,j=1,\ldots
\kappa\right\}.
\]
Because of the application we have in mind in the
next section, we will assume:
\begin{hypothesis}\ \\
a) $n$ is even and we set $n = 2p$\\[4pt]
b) The polynomials $p_j$ are of even degree and
their highest degree homogeneous parts have only
the origin as common zero.
\label{H}
\end{hypothesis}
One denotes by $C_+(K)$ the cone of polynomials
positive on $K$. In \cite[Theorem 1.4, p.
972]{putinar_indiana} it is proved that, under
Hypothesis \ref{H}, positive polynomials on $K$
belong to the additive cone
\[
C = \Sigma^2+p_1\Sigma^2+\cdots+p_\kappa\Sigma^2,
\]
where $\Sigma^2$ denotes the convex cone
generated by all squares of polynomials in
${\mathbb C}[x]$. The key result of
\cite{putinar_indiana} is:
\begin{theorem}(\cite[Lemma 3.2 p.
978]{putinar_indiana}). A functional $L$ on
${\mathbb R}[x]$ which is positive on $C$ is of
the form
\[
L(P)=\int_{K}P(x)d\nu(x),\quad P\in{\mathbb
R}[x],
\]
where $d\nu$ is a positive measure on $K$.
\end{theorem}

This result gives the solution to the moment problem on $K$: Let
$\alpha_{\ell,m}$ with $\ell, m\in{\mathbb N}_0^p$ be complex
numbers. Then there exists a positive measure on $K$ such that
\[
\int_K z^\ell z^{*m}d\nu(x)=\alpha_{\ell,m}
\]
if and only if the following conditions hold:
\begin{equation}
\begin{split}
L(|P(z,z^*)|^2)&\ge 0,\qquad \forall P\in{\mathbb
C}[x],\\
L(p_j(x) |s(z)|^2)&\ge 0,\qquad\forall s\in{\mathbb
C}[z],\quad j=1,\ldots p,\\
L(p_j(x)|P(z,z^*)|^2)&\ge 0, \quad j=p+1,\ldots
\kappa.
\end{split}
\end{equation}
\section{The case of a finite number of
generators}
\label{pg}
We now assume that the Abelian group
$\Gamma$ has a finite number, say $p$, of
generators, which we will denote by
$\gamma_1,\ldots, \gamma_p$. We assume that
there are independent in the sense that if
\[
\gamma_1^{n_1}\circ \cdots \circ\gamma_p^{n_p}= \iota
\]
for some integers $n_1,\ldots n_p\in{\mathbb Z}$,
then $n_1=\cdots=n_p=0$. In particular, each generator is of the form
\eqref{gamma},
$$\gamma_{i}(z)=\gamma_{\{\alpha_i\}}(z) =
(G_{\theta}\circ S_{\alpha_i}\circ G_{\theta}^{-1})(z)$$
with $\theta$ fixed, and where the set $\{\alpha_i\}_{i=1}^p$
generates a free discrete subgroup of the multiplicative group
of positive real numbers. We use in a free way
the multi-index notation.
\begin{theorem}
There is a positive measure $d\nu$ on the
distinguished boundary of the poly-disc such that
\[
\int_{\widehat{\Gamma}}\sigma(\gamma_1^{n_1})\cdots
\sigma(\gamma_p^{n_p})d\widehat{\mu}(\sigma)=
\int_{{\mathbb T}^p}e^{in_1\theta_1}\cdots
e^{in_p\theta_p}d\nu(\theta_1,\ldots, \theta_p).
\]
\label{nd}
\end{theorem}

To prove this theorem we specialize the results
of the preceding section to the case of the
poly-disc ${\mathbb D}^p$. It is a compact
algebraic set, with $n=2p$, $m=2p$ and
polynomials
\[
P_1(x)=1-|z_1|^2,\ldots, P_p(x)=1-|z_p|^2,
\]
and
\[
P_{p+1}(x)=|z_1|^2-1,\ldots, P_{2p}(x)=|z_p|^2-1.
\]

{\bf Proof of Theorem \ref{nd}:} We define a
linear form on polynomials in the variables
$z_1,\ldots, z_p, z_1^*,\ldots, z_p^*$ by
\[
L(z^\alpha z^{*\beta})=\int_{\widehat{\Gamma}}
(\sigma(\gamma_1))^{\alpha_1-\beta_1}
\cdots(\sigma(\gamma_p))^{\alpha_p-\beta_p}d\widehat{\mu}(\sigma).
\]
Let $p$ be a polynomial in the variables
$z_1,\ldots, z_p,
z_1^*,\ldots, z_p^*$. We write for short
\[
p(z,z^*)=p(z_1,\ldots, z_p, z_1^*,\ldots, z_p^*).
\]

Let $p(z,z^*)=\sum c_{\alpha,\beta}z^\alpha
z^{*{\beta}}$. Then
\[
\begin{split}
L(p(z,z^*))&=\\
&\hspace{-5mm}=\sum c_{\alpha,\beta}
\int_{\widehat{\Gamma}}(\sigma(\gamma_1)
)^{\alpha_1-\beta_1} \cdots (\sigma(\gamma_p)
)^{\alpha_p-\beta_p} d\widehat{\mu}(\sigma),
\end{split}
\]
and therefore we have
\begin{equation}
L(p(z,z^*))= \int_{\widehat{\Gamma}}p(
\sigma(\gamma_1),\sigma(\gamma_2),\ldots
,\sigma(\gamma_1)^*,\sigma(\gamma_2)^*,\ldots)
d\widehat{\mu}(\sigma).
\label{central}
\end{equation}

Since $|p(z,z^*)|^2$ is still a polynomial in $z$ and $z^*$, the
following conditions hold:
\begin{equation}
\label{mihai}
\begin{split}
L\left((1-z_jz_j^*)p(z,z^*)\right)&=0,
\quad j=1,2,\ldots p,\\
L\left(|p(z,z^*)|^2\right)&\ge 0.
\end{split}
\end{equation}
\mbox{}\qed\mbox{}\\

\begin{remark}
The fact that the characters are of modulus $1$
allows to prove \eqref{central} and
\eqref{mihai}. It does not seem possible to
relate our problems with another moment problem
when $p>1$ (for instance on the ball of
${\mathbb C}^p$).
\end{remark}
\begin{definition}
The Hermite transform of the element
\[
f(\sigma)=\sum_\alpha h_\alpha
\sigma(\gamma^\alpha)\] is
\[
{\mathbf I}(f)(z)=\sum_{\alpha} h_\alpha z^\alpha.
\]
\end{definition}
\begin{theorem}
Let $f\in{\mathbf L}_2(d\widehat{\mu})$ be such that the operator
of multiplication by $f$ defines a bounded operator from
${\mathbf L}_2(d\widehat{\mu})$ into itself. Then, for every
$g\in{\mathbf L}_2(d\widehat{\mu})$:
\begin{equation}
\label{ankara} {\mathbf I}(f g)={\mathbf I}(f){\mathbf I}(g).
\end{equation}
\end{theorem}
The proof is the same as for $p=1$ (see the proof of Theorem
\ref{hermite}). As in Definition \ref{dn:gtf}, the function of
$p+1$ variables
\[
{\mathscr H}(z,z_1,\ldots, z_p)=\sum_{n=0}^\infty z^n{\mathbf
I}(\widehat{h_n})(z_1,\ldots, z_p)
\]
is called the {\sl generalized transfer function} of the system.
\section{BIBO stability}
\label{bibo}
The system \ref{system} will be called bounded
input bounded output (BIBO) if there is an $M>0$
such that for every $\{u_n(\gamma)\}$ such that
\begin{equation}
\label{sup2} \sup_{n\in{\mathbb
N}_0}\|u_n(\cdot)\|_{\ell_2(\Gamma)}<\infty
\end{equation}
the output is such that
$\{y_n(\gamma)\}_{\gamma \in \Gamma}
\in \ell_2(\Gamma)$, $n=0,1,\ldots$, and it holds that
\begin{equation}
\label{BIBO} \sup_{n\in{\mathbb
N}_0}\|y_n(\cdot)\|_{\ell_2(\Gamma)}\le M \sup_{n\in{\mathbb
N}_0}\|u_n(\cdot)\|_{\ell_2(\Gamma)}.
\end{equation}
The following theorem gives a characterization of
BIBO systems. The proof follows the proof of
\cite[Theorem 3.2]{al_acap}. We note the
following difference between the two theorems:
in \cite{al_acap} the multiplication operators,
that is the counterparts of the operators
$\mathcal{M}_{h_n}$ defined here using the Wick
product, are automatically bounded. As explained
there, this is due to  V\r{a}ge's inequality
(see \cite[Proposition 3.3.2 p. 118]{MR1408433}
and (3.1) in \cite{al_acap}, and Section
\ref{table} below). Here we do not have an
analogue of this inequality.

\begin{theorem}
The system \eqref{system} is bounded input
bounded output if and only if the following two
conditions hold:\\
$(a)$ The multiplication operators \eqref{T_h}
\[
\mathcal{M}_{h_n}:\quad u\mapsto h_n\star u,\quad
n=0,1,\ldots
\]
are bounded from $\ell_2(\Gamma)$ into
itself.\\$(b)$ For all $v(\cdot)\in
\ell_2(\Gamma_+)$ with
$\|v(\cdot)\|_{\ell_2(\Gamma)}=1$ it holds that
\begin{equation}
\label{BIBO2} \sum_{n=0}^\infty
\|\mathcal{M}_{h_n}^*(v)\|_{\ell_2(\Gamma)}\le M.
\end{equation}
\label{istambul}
\end{theorem}

{\bf Proof:} That the condition \eqref{BIBO2} is
sufficient is readily seen. Indeed, take $v\in
\ell_2(\Gamma)$ with $\|v(\cdot)\|_{\ell_2(\Gamma)}=1$.
From \eqref{system} we have:
\begin{equation}
\label{y_nv}
\begin{split}
\langle y_n\, ,\, v\rangle_{\ell_2(\Gamma)}&=\sum_{m=0}^n\langle u_m\, ,\,
\mathcal{M}_{h_{n-m}}^*v\rangle_{\ell_2(\Gamma)},\quad n=0,1,\ldots,
\end{split}
\end{equation}
and hence
\[
\begin{split}
|\langle y_n\, ,\, v\rangle_{\ell_2(\Gamma)}| &\le
\sum_{m=0}^n \|u_m(\cdot) \|_{\ell_2(\Gamma)}\|
\mathcal{M}_{h_{n-m}}^*v\|_{\ell_2(\Gamma)}\\
&\le \left(\sup_{m=0,\ldots n}
 \|u_m(\cdot)\|_{\ell_2(\Gamma)}\right)\left(\sum_{m=0}^n
\|\mathcal{M}_{h_{n-m}}^*v\|_{\ell_2(\Gamma)}\right)
\\
&\le M\sup_{m\in{\mathbb N}_0}
 \|u_m(\cdot) \|_{\ell_2(\Gamma)}.
\end{split}
\]
We obtain \eqref{BIBO2} by taking
$v=y_n/\|y_n\|_{\ell_2(\Gamma)}$ when
$y_n\not=0$.\\

We now show that \eqref{BIBO2} is necessary. We
assume that the system is bounded input and
bounded output. We first note that the multiplication
operators \(\mathcal{M}_{h_n}\) are necessarily bounded.
Indeed, assume that \eqref{BIBO} is in force and
take $u_0=u\in \ell_2(\Gamma)$ and $u_n=0$
for $n > 0$. Then,
\[
y_n=h_n\star u=\mathcal{M}_{h_n}(u),\quad n=0,1,\ldots,
\]
and it follows from \eqref{BIBO} that
$\|\mathcal{M}_{h_n}\|\le M$ for $n=0,1,\ldots$\\

Let us now consider an input sequence $(u_n)$ which
satisfies \eqref{sup}. For a given $n$ and $v$
choose
\[
u_m=0\quad{\rm if}\quad \mathcal{M}_{h_{n-m}}^*v=0,
\]
and
\[
u_m=\frac{\mathcal{M}_{h_{n-m}}^*v}
{\|\mathcal{M}_{h_{n-m}}^*v\|_{\ell_2(\Gamma)}}\quad
{\rm otherwise}.
\]
We obtain from \eqref{y_nv} and \eqref{BIBO} that
\[
\sum_{m=0}^n\|\mathcal{M}_{h_{n-m}}^*v\|_{\ell_2(\Gamma)}\le M,
\]
from which we get \eqref{BIBO2}.
\mbox{}\qed\mbox{}\\

We now make a number of remarks: first,
condition \eqref{BIBO} is implied by the
stronger, but easier to deal with, condition
\begin{equation}
\label{BIBO3}
\sum_{n=0}^\infty \|{\mathcal
M}_{h_n}\|\le M.
\end{equation}
When $\Gamma$ is the trivial subgroup of
$SU(1,1)$, conditions \eqref{BIBO2} or
\eqref{BIBO3} reduce to the classical condition
\[
\sum_{n=0}^\infty |h_n|<\infty.
\]
Finally, other versions of this theorem could be
given, with non causal systems with respect to
the variable $n$ (as in \cite{al_acap}), or with
scale-causal signals. We state the last one. The proof
is the same as the proof of Theorem \ref{istambul}.

\begin{theorem}
The system \eqref{system} is scale-causal and
bounded input bounded output if and only if the
following two
conditions hold:\\
$(a)$ The multiplication operators \eqref{T_h}
\[
\mathcal{M}_{h_n}:\quad u\mapsto h_n\star u,\quad
n=0,1,\ldots
\]
are bounded from $\ell_2(\Gamma_+)$ into
itself.
\\ $(b)$ For all $v(\cdot)\in
\ell_2(\Gamma)$ with
$\|v(\cdot)\|_{\ell_2(\Gamma_+)}=1$ it holds that
\begin{equation*}
 \sum_{n=0}^\infty
\|\mathcal{M}_{h_n}^*(v)\|_{\ell_2(\Gamma_+)}\le
M.
\end{equation*}
\end{theorem}

\section{Dissipative systems}
\label{dissip}
We will call the system \eqref{system} {\sl
dissipative} if for every input sequence $(u_n)$
such that
\[
\sum_{n=0}^\infty
\|u_n(\cdot)\|_{{\ell_2(\Gamma)}}^2 <\infty
\]
it holds that
\begin{equation}
\label{bat-yam} \sum_{n=0}^\infty
\|y_n(\cdot)\|_{{\ell_2(\Gamma)}}^2\le \sum_{n=0}^\infty
\|u_n(\cdot)\|_{{\ell_2(\Gamma)}}^2.
\end{equation}
\begin{theorem}
The system is dissipative if and only if the
${\mathbf L}({\ell_2(\Gamma)})$-valued function
\[
S(z)=\sum_{n=0}^\infty z^n \mathcal{M}_{h_n}
\]
is analytic and contractive in the open unit
disc.
\end{theorem}

{\bf Proof:} Equations \eqref{bat-yam} expresses
that the block Toeplitz operator
\[
\begin{pmatrix}M_{h_0}&0&0&\cdots\\
               M_{h_1}&M_{h_0}&0&\cdots\\
               \vdots&\vdots& &\\
               & & &
               \end{pmatrix}
               \]
is a contraction from $\ell_2({\ell_2(\Gamma)})$ into itself, and
this is equivalent to the asserted condition on
$S$.\mbox{}\qed\mbox{}\\

We consider the case of
scale-causal signals (see Definition
\ref{scale_causal}).
\begin{definition}
The system \eqref{system} will be called
scale-causal dissipative if the following
conditions hold:\\
$(1)$ The operators $M_{h_n}$ are bounded from
$\ell_2(\Gamma_+)$ into itself.\\
$(2)$ Condition \eqref{bat-yam} holds, with
$\ell_2(\Gamma)$ replaced by $\ell_2(\Gamma_+)$.
\end{definition}

Recall that we have denoted
by ${\mathbf H}_2(d\nu)$ the closure in
${\mathbf L}_2(d\nu)$ of the powers $z^\alpha$,
where all the components of $\alpha$ are greater
or equal to $0$. Taking the Fourier and Hermite
transforms we have:
\begin{theorem}
\label{tm:scd}
The system is scale-causal dissipative if and only
the function
\begin{equation}
\label{izmir} {\mathscr H}(z,z_1,\ldots, z_p)=\sum_{n=0}^\infty
z^n{\mathbf I} ({\widehat{h_n}})(z_1,\ldots,z_p)
\end{equation}
is contractive from ${\mathbf H}_2({\mathbb
D})\otimes {\mathbf H}_2(d\nu)$ into itself.
Furthermore, if the space ${\mathbf H}_2(d\nu)$
is a reproducing kernel Hilbert space, say with
reproducing kernel $K(z_1,\ldots, w_1,\ldots)$,
condition \eqref{izmir} is equivalent to the
positivity of the kernel
\begin{equation}
\label{holon}
\frac{1-{\mathscr H} (z,z_1,\ldots){\mathscr H}
(w,w_1,\ldots)^*}{1-zw^*} K(z_1,\ldots, w_1,\ldots)
\end{equation}
in ${\mathbb D}^{p+1}$.
\end{theorem}

{\bf Proof:} Since the operators $M_{h_n}$ are assumed bounded,
we have
\[
h_{n-m}\star u_m\in\ell_2(\Gamma_+),\quad m=0,\ldots, n,
\]
for all entries $u_m\in \ell_2(\Gamma_+)$. Thus
\[
\widehat{h_{n-m}}\widehat{u_m}\in{\mathbf H}_2(d\widehat{\mu}),
\]
and we may apply Theorem \ref{hermite}. We can write:
\[
\left\|\sum_{n=0}^mh_{n-m}\star
u_m\right\|_{\ell_2(\Gamma_+)}=\left\|\sum_{n=0}^m{\mathbf
I}(\widehat{h_{n-m}}) {\mathbf
I}(\widehat{u_{m}})\right\|_{{\mathbf H}_2(d\nu)}.
\]
Thus the dissipativity is translated into the
contractivity of the block Toeplitz operator
\[
\begin{pmatrix}M_{\widehat{h_0}}&0&0&\cdots\\
M_{\widehat{h_1}}&M_{\widehat{h_0}}&0&\cdots\\
               \vdots&\vdots& &\\
               & & &
               \end{pmatrix}
               \]
from $\ell_2({\mathbf H}_2(d\nu))$ into itself, and hence the
claim on ${\mathscr H}$. To prove the second claim, we remark that
${\mathbf H}_2({\mathbb D})\otimes {\mathbf H}_2(d\nu)$ is the
reproducing kernel Hilbert space with reproducing kernel
\[
\frac{1}{1-zw^*} K(z_1,\cdots, w_1,\cdots).
\]
This comes from the fact that the reproducing
kernel of a tensor product of reproducing kernel
Hilbert spaces is the product of the reproducing
kernels; see \cite{aron}, \cite{saitoh}.
Condition \eqref{holon} follows then from the
well-known characterization of bounded
multipliers in reproducing kernel Hilbert
spaces; see for instance \cite{MR1839648}, and
the references therein.
\mbox{}\qed\mbox{}\\

\section{$\ell_1$-$\ell_2$ bounded systems}
\label{l1l2}
The system \eqref{system} will be
called $\ell_1$-$\ell_2$ bounded if there is a
$M>0$ such that for all inputs $(u_n)$ satisfying
\[
\sum_{n=0}^\infty\|u_n(\cdot)\|_{\ell_2(\Gamma)}<\infty,
\]
we have
\[
\left(\sum_{n=0}^\infty\|y_n(\cdot)
\|^2_{\ell_2(\Gamma)}\right)^{1/2} \le
M\sum_{n=0}^\infty\|u_n(\cdot)\|_{\ell_2(\Gamma)}.
\]
Taking the Fourier transform, this condition can
be rewritten as:
\begin{equation}
\label{l1l23}
\left(\sum_{n=0}^\infty\|\widehat{y_n}
\|^2_{{\mathbf L}_2(d\widehat{\mu})}\right)^{1/2}
\le M\sum_{n=0}^\infty\|\widehat{u_n}\|_{{\mathbf
L}_2(d\widehat{\mu})},
\end{equation}
The system \eqref{system} will be called {\sl
scale-causal $\ell_1$-$\ell_2$ bounded} if it is
moreover scale-causal, that is, if the operators
$M_{h_n}$ are bounded from $\ell_2(\Gamma_+)$
into itself. Condition \eqref{l1l23} then
becomes:
\begin{equation}
\label{l1l21}
\left(\sum_{n=0}^\infty\|\widehat{y_n}
\|^2_{{\mathbf H}_2(d\widehat{\mu})}\right)^{1/2}
\le M\sum_{n=0}^\infty\|\widehat{u_n}\|_{{\mathbf
H}_2(d\widehat{\mu})},
\end{equation}
from which we obtain, much in the same way as in
\cite{al_acap}, the following result. For
completeness we present a proof.
\begin{theorem}
A necessary and sufficient condition for the
system \eqref{system} to be scalar-causal and
$\ell_1$-$\ell_2$ bounded is that the function
\begin{equation}
\label{la_seine} H(z,\sigma)= \sum_{n=0}^\infty
z^n\widehat{h_n}(\sigma)\in{\mathbf H}_2({\mathbb
D})\otimes{\mathbf H}_2(d\widehat{\mu}),
\end{equation}
or, equivalently, that the transfer function
\begin{equation}
\label{eq:equiv} \mathscr H(z, z_1)= \sum_{n=0}^\infty
z^n{\mathbf I}( \widehat{h_n})(z_1)\in{\mathbf H}_2({\mathbb
D})\otimes{\mathbf H}_2(d\nu).
\end{equation}
\end{theorem}

{\bf Proof:} To see that condition
\eqref{la_seine} is necessary, it suffices to
take the sequence
\[
\widehat{u_n}=\begin{cases}1\quad{\rm if}\quad n=0\\
                 0\quad{\rm if}\quad n\not=0.
                 \end{cases}
\]
Then
\[
\widehat{y_n}=\widehat{h_n}\quad n=0,1,\ldots,
\]
and condition \eqref{l1l21} implies that
$H(z,\sigma)\in{\mathbf H}_2({\mathbb
D})\otimes{\mathbf H}_2(d\widehat{\mu})$.
Conversely, assume that the function
$H(z,\sigma)\in{\mathbf H}_2({\mathbb
D})\otimes{\mathbf H}_2(d\widehat{\mu})$. From
the expression \eqref{conv}, and using the
Cauchy-Schwarz inequality on
\[
Y(z,\sigma)=H(z,\sigma)U(z,\sigma)
=\sum_{n=0}^\infty
H(z,\sigma)\left(z^n\widehat{u_n}\right),
\]
we have
\[
\begin{split}
\|Y(z,\sigma)\|_{{\mathbf H}_2({\mathbb
D})\otimes{\mathbf H}_2(d\widehat{\mu})}&\le
\sum_{n=0}^\infty\| H(z,\sigma)z^n\widehat{u_n}
\|_{{\mathbf
H}_2({\mathbb
D})\otimes{\mathbf H}_2(d\widehat{\mu})}\\
&\le \sum_{n=0}^\infty\|z^n\widehat{u_n}
\|_{{\mathbf H}_2({\mathbb
D})\otimes{\mathbf
H}_2(d\widehat{\mu})} \|H(z,\sigma)|
\|_{{\mathbf H}_2({\mathbb
D})\otimes{\mathbf
H}_2(d\widehat{\mu})}.
\end{split}
\]
But we have that
\[
\|z^n\widehat{u_n} \|_{{\mathbf H}_2({\mathbb
D})\otimes{\mathbf H}_2(d\widehat{\mu})}=
\|\widehat{u_n} \|_{{\mathbf
H}_2(d\widehat{\mu})},
\]
and so we obtain \eqref{l1l21} with $M=
 \|H(z,\sigma)|
\|_{{\mathbf H}_2({\mathbb D})\otimes{\mathbf
H}_2(d\widehat{\mu})}$. The equivalence with condition
\eqref{eq:equiv} follows by taking the Hermite transform.
\mbox{}\qed\mbox{}\\

When ${\mathbf H}_2(d\nu)\not ={\mathbf L}_2(d\nu)$ (recall that
$\Gamma$ is finitely generated), \eqref{eq:equiv} can be
translated into reproducing kernel conditions. In particular, in
the cyclic case, we have:
\begin{theorem}
Assume that ${\mathbf H}_2(d\nu)\not = {\mathbf L}_2(d\nu)$, and
let
\[
\frac{A(z_1)A(w_1)^*-B(z_1)B(w_1)^*}{1-z_1w_1^*}
\]
be the reproducing kernel of ${\mathbf
H}_2(d\nu)$.
 The system
\eqref{system} is scale-causal and
$\ell_1-\ell_2$ bounded if and only if there is
a $M>0$ such that the kernel
\[
\frac{A(z_1)A(w_1)^*-B(z_1)B(w_1)^*}{(1-zw^*)(1-z_1w_1^*)}-
M{\mathscr
  H}(z,z_1)
{\mathscr H}(w,w_1)^*
\]
is positive in the bi-disc.
\end{theorem}

As in the case of equation \eqref{holon}, this
comes from the characterization of the
reproducing kernel of a tensor product of
reproducing kernel Hilbert spaces.

\section{The white noise space setting and a table}
\label{table}
Another kind of double convolution system, with a setting quite
similar to the setting presented here, has been developed in
\cite{al_acap}, and rely on Hida's theory of the white noise space
(see \cite{MR1408433}, \cite{hida_taiwan}, \cite{MR1244577} for
the latter). We now review the main features of Hida' theory and
of \cite{al_acap}. The starting point in Hida's theory is the
function
\[
K(s_1-s_2)=e^{-{\frac{\|s_1-s_2\|_{{\mathbf
L}_2({\mathbb R})}^2}{2}}},
\]
which is positive in the sense of reproducing kernels for
$s_1,s_2$ in the Schwartz space ${\mathcal S}$ of real valued
rapidly vanishing functions. By the Bochner-Minlos theorem there
exists a probability measure $P$ on the dual space ${\mathcal
S}^\prime$ of real valued tempered distributions such that
\[
e^{-\frac{{\|s\|_{{\mathbf L}_2({\mathbb
R})}^2}}{2}}=\int_{{\mathcal
S}^\prime}e^{i\langle w,s\rangle} dP(w) ,\quad
s\in{\mathcal S}^\prime,
\]
where we have denoted by $\langle w,s\rangle$
the duality between ${\mathcal S}$ and
${\mathcal S}^\prime$. The white noise space is
defined to be the real Hilbert space ${\mathbf
L}_2({\mathcal S}^\prime,{\mathcal F},P)$, where
${\mathcal F}$
denotes the underlying Borelian sigma-algebra.\\

Among all orthonormal basis of ${\mathbf L}_2({\mathcal
S}^\prime,{\mathcal F},P)$, there is one which plays a special
role; it is constructed from the Hermite functions, and is indexed
by the set $\ell$ of infinite sequences
$(\alpha_1,\alpha_2,\ldots)$ indexed by ${\mathbb N}$, and with
values in ${\mathbb N}_0$, and for which $\alpha_j=0$ for all $j$
at the exception of at most a finite number of $j$. See
\cite[Definition 2.2.1 p. 19]{MR1408433}. We will denote by
$H_\alpha$ (with $\alpha\in\ell$) the elements of this basis. An
element
\begin{equation}
\label{thessa} F=\sum_{\alpha\in\ell}f_\alpha
H_\alpha,\quad f_\alpha\in{\mathbb R},
\end{equation}
belongs to the white noise space if
\[
\sum_{\alpha\in\ell}\alpha!f_\alpha^2<\infty.
\]
The Wick product is defined by
\[
H_\alpha\lozenge H_\beta=H_{\alpha+\beta},\quad
\alpha,\beta\in\ell.
\]
The white noise space is not stable under the
Wick product, and there is the need to introduce
a nuclear space, called the Kondratiev space,
within which the Wick product is stable. The
Kondratiev space is the projective limit of the
real Hilbert spaces ${\mathcal H}_k$ of formal
sums of the form \eqref{thessa} for which
\[
\sum_\alpha f^2_\alpha(2{\mathbb
N})^{-q\alpha}<\infty \label{Kodd}
\]
for some $q\in{\mathbb N}$, where we use the
notation
\[
(2{\mathbb N})^\alpha\stackrel{\rm def.}{=}2^{\alpha_1}\times
4^{\alpha_2}\times 6^{\alpha_3}\cdots
\]
One can also consider the complexified versions
of these spaces.\\

We also recall
 V\r{a}ge's inequality
(see \cite[Proposition 3.3.2 p. 118]{MR1408433}):
Fix some integer $l>0$,  and let $k>l+1$.
Consider $h\in {\mathcal H}_{l}$ and $u\in
{\mathcal H}_{k}$. Then, $h\lozenge
u\in{\mathcal H}_k$ and
\begin{equation*}
\|h\lozenge u\|_{k}\le
A(k-l)\|h\|_{l}\|u\|_{k},
\end{equation*}
where
\begin{equation*}
A(k-l)=\sum_{\alpha\in\ell}(2{\mathbb
N})^{(l-k)\alpha}
\end{equation*}
is a finite number.\\

We can now introduce the systems considered in \cite{al_acap}. A
system will be characterized by a sequence $(h_n)_{n\in{\mathbb
Z}}$ of elements in ${\mathcal H}_l$ for some $l\in{\mathbb N}$,
and a signal will be a sequence of elements in one of the spaces
${\mathcal H}_k$, with $k>l+1$. Input-output relations are
expressions of the form
\[
y_n=\sum_{m\in{\mathbb Z}}h_{n-m}\lozenge x_m,\quad n\in\mathbb Z.
\]
Note that in view of V\r{a}ge's inequality the
output sequence consists also of elements of
${\mathcal H}_k$. Furthermore, decomposing this
equation along the basis $H_\alpha$ we obtain
the double convolution system
\begin{equation*}
y_\alpha(n)=\sum_{m\in{\mathbb Z}}\sum_{\beta\le \alpha}
h_{\alpha-\beta}(n-m)u_\beta(m),\quad n\in{\mathbb Z}.
\end{equation*}
The map ${\mathbf I}$ which to $H_\alpha$
associates the polynomial $z^\alpha$ is called
the Hermite transform. It is such that
\[
{\mathbf I}(f\lozenge g)={\mathbf I}(f) {\mathbf
I}(g),\quad \forall f,g\in S_{-1}.
\]
Note that under the Hermite transform the white
noise space is mapped onto the reproducing kernel
Hilbert space with reproducing kernel
$e^{\langle z,w\rangle_{\ell_2}}$, that is, onto
the Fock space.\\

\newpage
We now give the table presenting the parallels between the white
noise space case (as applied in the paper \cite{al_acap}), and
the present multi-scale case. The reader might want to look at a
similar table in \cite{al_pota2}, where the analogies between the
white noise space case and the hyper-holomorphic case are
presented.

\begin{tabular}{|m{0.2\linewidth}|>{\centering}m{0.35\linewidth}|m{0.35\linewidth}|}
\hline & &\\
{The setting}&{Stochastic case}&$\qquad$ {Multi-scale case}
\\
 & &\\\hline
& &\\
Underlying space & The white noise space &
$\qquad\qquad \ell_2(\Gamma)$\\
& &\\ \hline& &\\
Hermite transform& $\mathbf I(H_\alpha)=z^\alpha$&
{\centering{$\mathbf I(\sigma(\gamma)^\alpha)=z^\alpha$\\ ($\Gamma$: finitely generated)}}
\\
& &\\
\hline& &\\
 Underlying reproducing kernel Hilbert space&The Fock space
 & $\quad$ The space ${\mathbf H}_2(d\nu)$\\
  &&\\
 \hline& &\\
 Key tool used&{Minlos  theorem \\(to build the white noise space)}
 &Moment problem on the
poly-disc (to build the Hermite
transform) \\
 & &\\
 \hline& &\\
The product& Wick product&{\centering{Convolution\\ $\quad$ with
respect to $\Gamma$.}}\\
 & &\\
\hline & &\\
 {Double convolution}  &
{\hspace*{-7mm}$y_\alpha(n)=\sum_{m\in{\mathbb
Z}}\sum_{\beta\le \alpha}$\\$\qquad h_{\alpha-\beta}(n-m)u_\beta(m)
$} &{\centering{\hspace*{-9mm} $y_n(\gamma)=\sum_{m=0}^n
\sum_{\varphi\in\Gamma}$\\$\qquad h_{n-m}(\gamma\circ\varphi^{-1})u_m(\varphi)$}}
\\
& &\\
 \hline
\end{tabular}
\mbox{}\\

\begin{remark}
The pointwise product in ${\mathbf
L}_2(\widehat{d\mu})$ is a convolution in
$\ell_2(\Gamma)$. Strictly speaking, it would be
better to define the Hermite transform as the
composition of the Fourier transform and of the
map $H_\alpha\mapsto z^\alpha$.
\end{remark}

\end{document}